\documentclass[a4paper,12pt]{amsart}
\usepackage{amsfonts}
\usepackage{amsmath}
\usepackage{amssymb}
\usepackage[a4paper]{geometry}
\usepackage{mathrsfs}
\usepackage{hyperref}
\renewcommand\eqref[1]{(\ref{#1})} 
\numberwithin{equation}{section}
\theoremstyle{plain}
\newtheorem{theorem}{Theorem}[section]

\theoremstyle{definition}
\newtheorem{definition}[theorem]{Definition}
\newtheorem{remark}[theorem]{Remark}

\begin{document}
\title[Blow-up results for the NLS equations]
{Blow-up results for the semilinear Schrödinger equations with forcing and gradient terms: the critical cases}

\author[B. T. Torebek]{Berikbol T. Torebek}

\address{Berikbol T. Torebek \newline Institute of
Mathematics and Mathematical Modeling,  Almaty, Kazakhstan \newline and \newline Department of Mathematics: Analysis, Logic and Discrete Mathematics \newline Ghent University, Ghent, Belgium}
\email{{berikbol.torebek@ugent.be}}



\begin{abstract}
The paper is devoted to the study of critical cases of the nonlinear Schrödinger (NLS) equation with source and gradient terms, subsequently providing answers to some open questions posed by Alotaibi et al in [Z. Angew. Math. Phys., 73 (2022), 1--17]. The main results state that in critical cases, the problem under consideration does not have any global-in-time weak solutions. The proof approach is based on
modified method of test functions specifically adapted to the nature of considered problems.
\end{abstract}

\maketitle

\section{Introduction and problem formulation}
This paper investigates the Cauchy problem to the semilinear Schr\"{o}dinger equation
\begin{equation}
\label{1}
\begin{cases}
iu_t+\Delta_{x}u=\lambda|u|^{p}+\mu|\nabla u|^q+f(x),\,\,\,(0,\infty)\times\mathbb{R}^{n},\\
u(0, x)=u_{0}(x),\,\,\,\,\, x\in\mathbb{R}^{n},
\end{cases}
\end{equation}
where $u_0(x)$ and $f(x)$ are given complex-valued functions and $\lambda\in\mathbb{C}\setminus\{0\}, \mu\in\mathbb{C}$.

In the cases $f\equiv 0$ and $\mu=0$, i.e.
\begin{equation}\label{NLS}iu_t+\Delta_{x}u=\lambda|u|^{p},\,\,\,(0,\infty)\times\mathbb{R}^{n},\end{equation} the problem \eqref{1}
has been studied by many authors (see \cite{Tao, Caze, Fuj, Fuj1, Ikeda1, Ikeda2, Ikeda, Jleli, Kirane, Kirane1, Oh}). For example, in \cite{Ikeda} Ikeda and Wakasugi studied the critical exponent in the sense of Fujita for \eqref{NLS} and proved that when $1<p\leq 1+\frac{2}{n},$ problem \eqref{NLS} does not have any global $L^2$-solutions. Later, Ikeda and Inui proved the $L^2$-small data blow-up results for $1<p<1+\frac{4}{n}$ (see \cite{Ikeda1}) and $L^2$-large data blow-up results for $p>1+\frac{4}{n}$ (see \cite{Ikeda2}). Kirane et al. (see \cite{Kirane, Kirane1}) studied nonlocal extensions of problem \eqref{NLS}.

The problem \eqref{1} was studied in \cite{Alo} by Alotaibi et al. and the following results were obtained:
\begin{itemize}
    \item[(i)] Let $u_0,f\in L^1_{loc}(\mathbb{R}^n),$ $\lambda\in\mathbb{C}\setminus\{0\},$ $\mu=0,$ and $$p^*=\left\{\begin{array}{cc}
        \infty, & n=1,2, \\
        \frac{n}{n-2}, & n\geq 3.
    \end{array}\right.$$ If $\lambda_i\int_{\mathbb{R}^n}f_i(x)dx>0, \,i=1,2$, and $1<p<p^*,$ then problem \eqref{1} does not admit a weak global solution. However, if $p>p^*$, then \eqref{1} admits global solutions for some $u_0$ and $f.$ Here $\lambda_1=\Re(\lambda), \lambda_2=\Im(\lambda),$ $f_1(x)=\Re(f(x)), f_2(x)=\Im(f(x)).$
    \item[(ii)] Let $u_0, \nabla u_0,f\in L^1_{loc}(\mathbb{R}^n),$ $\lambda,\mu\in\mathbb{C}\setminus\{0\},$ $\mu_1=\Re(\mu),\,\mu_2=\Im(\mu)$ and $$q^*=\left\{\begin{array}{cc}
        \infty, & n=1, \\
        \frac{n}{n-1}, & n\geq 2.
    \end{array}\right.$$
    If $$1<p<p^*,\,\lambda_i\mu_i\geq 0,\,\lambda_i\int_{\mathbb{R}^n}f_i(x)dx>0, \,i=1,2$$ or $$1<q<q^*,\,\lambda_i\mu_i> 0,\,\mu_i\int_{\mathbb{R}^n}f_i(x)dx>0, \,i=1,2,$$ then \eqref{1} admits no global weak solution. If $n\geq 3,$ $p>p^*$ and $q>q^*,$ then \eqref{1} admits global solutions (stationary solutions) for some $u_0$ and $f.$
\end{itemize}
From (ii), one can deduce that the Fujita-type critical exponent for \eqref{1} is given by
$$p_{c}(n,q)=\left\{\begin{array}{cc}
        \infty, & 1<q<q^*, \\
        \frac{n}{n-2}, & q\geq q^*.
    \end{array}\right.$$
One observes that the gradient nonlinearity induces an interesting phenomenon of discontinuity of the critical exponent, jumping from
$\frac{n}{n-2}$ to $\infty$ as $q$ reaches the value from above $\frac{n}{n-1}$. A similar phenomenon for parabolic equations with combined nonlinearities was previously discovered in \cite{Souplet, Torebek}.

It should be noted that in (i) and (ii) there are no results for the critical cases $p=p^*,\,q>1$ and $q=q^*,\,p>1.$ Therefore, the main objective of this paper is to study critical cases of problem \eqref{1} which complement the results of (i) and (ii).
\subsection{Main results}
Before formulating the main result, we give the definition of a weak solution to problem \eqref{1}.
\begin{definition}\label{Def1}
A function $u$ is called a weak solution of the problem \eqref{1} with $\mu=0,$ if $u\in L^{p}_{loc}((0,T)\times\mathbb{R}^{n})$ and satisfies the identity
\begin{multline}\label{2}
i\int_{\mathbb{R}^{n}}u_{0}(x)\varphi(0, x)dx+\int\limits_0^T\int\limits_{\mathbb{R}^n}\left[\lambda|u(t, x)|^{p}+f(x)\right]\varphi(t, x)dxdt\\
=\int\limits_0^T\int\limits_{\mathbb{R}^n} u(t, x)\Delta_{x}\varphi(t, x)dxdt-i\int\limits_0^T\int\limits_{\mathbb{R}^n}u(t, x)\varphi_t(t, x)dxdt,
\end{multline}
for $u_0, f\in L^{1}_{loc}(\mathbb{R}^n)$ and for any compactly supported $0\leq \varphi\in C^{1, 2}_{t, x}((0,T)\times\mathbb{R}^{n})$.
\end{definition}
\begin{definition}\label{Def2}
A function $u$ is called a weak solution of the problem \eqref{1} if $u\in L^{p}_{loc}((0,T)\times\mathbb{R}^{n}), \nabla u\in L^{q}_{loc}((0,T)\times\mathbb{R}^{n})$ and satisfies the identity
\begin{multline}\label{2+}
i\int_{\mathbb{R}^{n}}u_{0}(x)\varphi(0, x)dx+\int\limits_0^T\int\limits_{\mathbb{R}^n}\left[\lambda|u(t, x)|^{p}+\mu|\nabla u(t,x)|^q+f(x)\right]\varphi(t, x)dxdt\\
=-\int\limits_0^T\int\limits_{\mathbb{R}^n} \nabla u(t, x)\nabla\varphi(t, x)dxdt-i\int\limits_0^T\int\limits_{\mathbb{R}^n}u(t, x)\varphi_t(t, x)dxdt,
\end{multline}
for $u_0, f\in L^{1}_{loc}(\mathbb{R}^n)$ and for any compactly supported $0\leq \varphi\in C^{1, 2}_{t, x}((0,T)\times\mathbb{R}^{n})$.
\end{definition}
Now, we formulate the main result.
\begin{theorem}\label{Th1}
Let $u_0, f\in L^{1}_{loc}(\mathbb{R}^n),$ $\lambda\in \mathbb{C}\setminus\{0\}$ and $u$ is the solution of problem \eqref{1}.
\begin{itemize}
    \item[(a)] If $\mu=0,$ and $$p=\frac{n}{n-2}\,\,\,\,\text{and}\,\,\,\,\lambda_i\int_{\mathbb{R}^n}f_i(x)dx>0, \,i=1,2,$$ then problem \eqref{1} does not admit a weak global solution.
    \item[(b)] Let $\mu\in \mathbb{C}\setminus\{0\}.$ If
    $$p=p^*=\frac{n}{n-2},\,q>1,\,\,\,\,\text{and}\,\,\,\,\lambda_i\mu_i\geq 0,\,\lambda_i\int_{\mathbb{R}^n}f_i(x)dx>0, \,i=1,2$$ or $$q=q^*=\frac{n}{n-1},\,p>1,\,\,\,\,\text{and}\,\,\,\,\lambda_i\mu_i> 0,\,\mu_i\int_{\mathbb{R}^n}f_i(x)dx>0, \,i=1,2,$$ then \eqref{1} admits no global weak solution.
\end{itemize}
\end{theorem}
\begin{remark}
Analyzing the results of Theorem \ref{Th1}, we can draw the following conclusions:
\begin{itemize}
    \item[(I)] Part (a) of Theorem \ref{Th1} provides an answer to the open question in \cite{Alo} (see Remark 1.2) and therefore complements Theorem 1.1 of that paper. Accordingly, part (b) of Theorem \ref{Th1} complements Theorem 1.3 in \cite{Alo}.
    \item[(II)] In Theorem \ref{Th1}, no additional properties $u_0$ are required other than $u_0\in L^1_{loc}(\mathbb{R}^n)$. This means that Theorem \ref{Th1} is true even for $u_0(x)=0,\,x\in\mathbb{R}^n$. In this connection, the following question arises: Is it possible to study the critical exponents of problem \eqref{1} if $$\int_{\mathbb{R}^n}f_i(x)dx=0,\,f_i(x)\not\equiv 0?$$ Unfortunately, this question remains open even for the case of parabolic equation, due to technical difficulties.
\end{itemize}
\end{remark}

\section{Proof of Theorem \ref{Th1}}
The proof will be carried out using the test function method. However, unlike the classical method of test functions (see, for example, \cite{Zhang, Zhang1, Pohoz}), we will use a modification based on test functions with logarithmic arguments. This approach has previously been used successfully in \cite{Bori, Samet, Torebek}.
\subsection{The proof of case (i)}
Note that the integral identity \eqref{2} is equivalent to the system of integral equations
\begin{equation}\label{def1}\begin{split}
&-\int_{\mathbb{R}^{n}}u_{02}(x)\varphi(0, x)dx+\int\limits_0^T\int\limits_{\mathbb{R}^n}\left[\lambda_1|u(t, x)|^{p}+f_1(x)\right]\varphi(t, x)dxdt\\&
=\int\limits_0^T\int\limits_{\mathbb{R}^n} u_1(t, x)\Delta\varphi(t, x)dxdt+\int\limits_0^T\int\limits_{\mathbb{R}^n}u_2(t, x)\varphi_t(t, x)dxdt,
\end{split}
\end{equation}
and
\begin{align*}
&\int_{\mathbb{R}^{n}}u_{01}(x)\varphi(0, x)dx+\int\limits_0^T\int\limits_{\mathbb{R}^n}\left[\lambda_2|u(t, x)|^{p}+f_2(x)\right]\varphi(t, x)dxdt\\&
=\int\limits_0^T\int\limits_{\mathbb{R}^n} u_2(t, x)\Delta\varphi(t, x)dxdt-\int\limits_0^T\int\limits_{\mathbb{R}^n}u_1(t, x)\varphi_t(t, x)dxdt,
\end{align*}
where $u_{01}=\Re(u_0)$ and $u_{02}=\Im(u_0).$

The proof will be given by contradiction. Suppose that there exists a non-trivial global solution. Let us consider the case $$\lambda_1\int_{\mathbb{R}^n}f_1(x)dx>0.$$ Then from \ref{def1} we have
\begin{equation}\label{3}\begin{split}
\int\limits_0^T\int\limits_{\mathbb{R}^n}\left[|u(t, x)|^{p}+\frac{1}{\lambda_1}f_1(x)\right]\varphi(t, x)dxdt&
\leq \frac{1}{|\lambda_1|}\int_{\mathbb{R}^{n}}|u_0(x)|\varphi(0, x)dx\\+\frac{1}{|\lambda_1|}\int\limits_0^T\int\limits_{\mathbb{R}^n} |u(t, x)||\Delta\varphi(t, x)|dxdt&+\frac{1}{|\lambda_1|}\int\limits_0^T\int\limits_{\mathbb{R}^n}|u(t, x)||\varphi_t(t, x)|dxdt,\end{split}
\end{equation}
We set the test function as $$\varphi(t,x)=\eta(t)\phi(x),$$
for large enough $R, T,$ where
\begin{equation}\label{T2}
\eta(t)=\nu^l\left(\frac{t}{T}\right),\,l>\frac{2p}{p-1},\,\,\, t>0, \end{equation}
and
\begin{equation}\label{T1}
\phi(x)=\Psi^m\left(\frac{\ln\left(\frac{|x|}{\sqrt{R}}\right)}{\ln\left(\sqrt{R}\right)}\right),\,m>\frac{2p}{p-1},\,\, x\in\mathbb{R}^{n}.
\end{equation}
Let $\nu\in C^1(\mathbb{R})$ be such that
$\nu\geq0;\,\,\,\nu\not\equiv0;\,\,\,\text{supp}(\nu)\subset(0,1),$ and $\Psi\in C^2(\mathbb{R}^n)$ be a smooth function  satisfying
\begin{equation}\label{TT}
\Psi(s)=\left\{\begin{array}{l}
1,\,\,\,\,\,\text{if}\,\,-\infty< s\leq0,\\
\searrow,\,\,\text{if}\,\,0< s<1,\\
0,\,\,\,\,\,\text{if}\,\,s\geq1.\end{array}\right.
\end{equation}
Using the fact that  $\text{supp}(\nu)\subset(0,1)$, we can easily get
\begin{equation*}
\int\limits_{\mathbb{R}^n}|u_{0}(x)|\varphi(0,x)dx=\nu(0)\int\limits_{\mathbb{R}^{n}}|u_{0}(x)|\phi(x)dx=0.
\end{equation*}
Also from property of $\eta$ it immediately follows that
\begin{equation}\label{f}
\int\limits_0^T\int\limits_{\mathbb{R}^n} f_1(x)\varphi dx dt=\int\limits_0^T\nu^l\left(\frac{t}{T}\right)dt\int\limits_{\mathbb{R}^n} f_1(x)\phi(x) dx =CT\int\limits_{\mathbb{R}^n} f_1(x)\phi(x) dx
\end{equation}
By using the $\varepsilon$-Young inequalities to the right hand side of \eqref{3} we have
\begin{align*}
\int\limits_0^T\int\limits_{\mathbb{R}^n}|u||\varphi_t|dxdt&\leq \varepsilon_1\int\limits_0^T\int\limits_{\mathbb{R}^n} |u|^p\varphi dxdt +C_1(\varepsilon)\int\limits_0^T\int\limits_{\mathbb{R}^n}|\varphi_t|^{\frac{p}{p-1}}{|\varphi|^{-\frac{1}{p-1}}}dxdt,
\end{align*}
and
\begin{align*}
\int\limits_0^T\int\limits_{\mathbb{R}^n}|u||\Delta\varphi| dxdt&\leq \varepsilon_2\int\limits_0^T\int\limits_{\mathbb{R}^n}|u|^p\varphi dxdt+C_2(\varepsilon_2)\int\limits_0^T\int\limits_{\mathbb{R}^n}{|\Delta \varphi|^{\frac{p}{p-1}}}{|\varphi|^{-\frac{1}{p-1}}}dxdt,
\end{align*} where $C_1$ and $C_2$ are some positive constants.

Choosing $\varepsilon_1=\varepsilon_2=\frac{|\lambda_1|}{2}$ and combining the above calculations with \eqref{3} and \eqref{f} one can get
\begin{equation}\label{T3}\begin{split}
&\frac{1}{\lambda_1}\int\limits_{\mathbb{R}^n} f_1(x)\phi(x) dx dt\leq CT^{-1}\left(I_1+I_2\right),
\end{split}\end{equation} with \begin{equation*}\begin{split}
&I_1=\int\limits_0^T\int\limits_{\mathbb{R}^n}|\varphi_t|^{\frac{p}{p-1}}{|\varphi|^{-\frac{1}{p-1}}}dxdt=\int_0^T\eta^{-\frac{1}{p-1}}\left|\eta_t\right|^{\frac{p}{p-1}}dt\int_{\mathbb{R}^{n}}\phi(x)dx,
\\&I_2=\int\limits_0^T\int\limits_{\mathbb{R}^n}{|\Delta \varphi|^{\frac{p}{p-1}}}{|\varphi|^{-\frac{1}{p-1}}}dxdt=\int_0^T\eta(t)dt\int_{\mathbb{R}^{n}}\phi^{-\frac{1}{p-1}}\left|\Delta\phi\right|^{\frac{p}{p-1}}dx.
\end{split}\end{equation*}
Using the properties of $\eta$ and $\phi,$ as well as changing $t=T\tau,\, x={\sqrt{R}} y,$ and applying $|\phi'(s)|\leq C,\,|\phi''(s)|\leq C$ we have that
\begin{equation*}\begin{split}
&I_1\leq C T^{1-\frac{p}{p-1}}R^n,\,\,\,\,\text{and}\,\,\,\,I_2\leq C T (\ln R)^{-\frac{p}{p-1}}.
\end{split}\end{equation*}
Since, $p=\frac{n}{n-2}$, by inserting the last inequalities into \eqref{T3} and choosing $T=R^4$, we can verify that
\begin{equation}\label{T4}\begin{split}
&\frac{1}{\lambda_1}\int\limits_{\mathbb{R}^n} f_1(x)\phi(x) dx\leq C\left(R^{-n}+(\ln R)^{-\frac{n}{2}}\right),
\end{split}\end{equation}
Hence passing to the limit as $R\to\infty$ in \eqref{T4} and taking into account $$\lim\limits_{R\rightarrow+\infty}\phi\left(x\right)=\lim\limits_{R\rightarrow+\infty}\Psi^m\left(\frac{\ln\left(\frac{|x|}{{\sqrt{R}}}\right)}{\ln\left({\sqrt{R}}\right)}\right)=\Phi(-1)=1,$$ we deduce that \begin{align*}
\frac{1}{\lambda_1}\lim\limits_{R\rightarrow +\infty}\int\limits_{\mathbb{R}^n} f_1(x)\phi(x) dx&=\frac{1}{\lambda_1}\int\limits_{\mathbb{R}^n} f_1(x) dx\leq 0,\end{align*} which is a contradiction with ${\lambda_1}\int\limits_{\mathbb{R}^n} f_1(x) dx>0$.

The case $\lambda_2\int_{\mathbb{R}^n}f_2(x)dx>0$ can be proved in a similar way.

\subsection{The proof of case (ii)} It is not difficult to show that \eqref{2+} is equivalent to the system of integral equations (see \cite[Remark 1.4]{Jleli})
\begin{equation}\label{def2}\begin{split}
&-\int_{\mathbb{R}^{n}}u_{02}(x)\varphi(0, x)dx+\int\limits_0^T\int\limits_{\mathbb{R}^n}\left[\lambda_1|u(t, x)|^{p}+\mu_1|\nabla u(t,x)|^q+f_1(x)\right]\varphi(t, x)dxdt\\&
=\int\limits_0^T\int\limits_{\mathbb{R}^n} u_1(t, x)\Delta\varphi(t, x)dxdt+\int\limits_0^T\int\limits_{\mathbb{R}^n}u_2(t, x)\varphi_t(t, x)dxdt,
\end{split}
\end{equation}
and
\begin{align*}
&\int_{\mathbb{R}^{n}}u_{01}(x)\varphi(0, x)dx+\int\limits_0^T\int\limits_{\mathbb{R}^n}\left[\lambda_2|u(t, x)|^{p}+\mu_2|\nabla u(t,x)|^q+f_2(x)\right]\varphi(t, x)dxdt\\&
=\int\limits_0^T\int\limits_{\mathbb{R}^n} u_2(t, x)\Delta\varphi(t, x)dxdt-\int\limits_0^T\int\limits_{\mathbb{R}^n}u_1(t, x)\varphi_t(t, x)dxdt.
\end{align*}
Suppose that $u$ is a global weak solution to problem \eqref{1}. 

Let us consider the case
$p=p^*=\frac{n}{n-2},\,q>1,$ and $\lambda_1\mu_1\geq0,\,\lambda_1\int_{\mathbb{R}^n}f_1(x)dx>0.$ Then from \ref{def2} we have
\begin{equation*}\begin{split}
&\int\limits_0^T\int\limits_{\mathbb{R}^n}\left[|u(t, x)|^{p}+\frac{\mu_1}{\lambda_1}|\nabla u(t,x)|^q+\frac{1}{\lambda_1}f_1(x)\right]\varphi(t, x)dxdt
\\&\leq \frac{1}{|\lambda_1|}\int\limits_0^T\int\limits_{\mathbb{R}^n} |u(t, x)||\Delta\varphi(t, x)|dxdt+\frac{1}{|\lambda_1|}\int\limits_0^T\int\limits_{\mathbb{R}^n}|u(t, x)||\varphi_t(t, x)|dxdt\\&+\frac{1}{|\lambda_1|}\int\limits_{\mathbb{R}^{n}}|u_0(x)|\varphi(0, x)dx,\end{split}
\end{equation*}
Since $\lambda_1\mu_1\geq0,$ the last inequality can be rewritten as
\begin{equation*}\begin{split}
\int\limits_0^T\int\limits_{\mathbb{R}^n}\left[|u(t, x)|^{p}+\frac{1}{\lambda_1}f_1(x)\right]\varphi(t, x)dxdt
&\leq \frac{1}{|\lambda_1|}\int\limits_{\mathbb{R}^{n}}|u_0(x)|\varphi(0, x)dx\\+\frac{1}{|\lambda_1|}\int\limits_0^T\int\limits_{\mathbb{R}^n} |u(t, x)||\Delta\varphi(t, x)|dxdt&+\frac{1}{|\lambda_1|}\int\limits_0^T\int\limits_{\mathbb{R}^n}|u(t, x)||\varphi_t(t, x)|dxdt.\end{split}
\end{equation*}
Then, completely repeating the calculations carried out in the proof of part (i), we obtain a contradiction with $\lambda_1\int_{\mathbb{R}^n}f_1(x)dx>0.$

The case $p=p^*=\frac{n}{n-2},\,q>1,$ $\lambda_2\mu_2\geq0,$ $\lambda_2\int_{\mathbb{R}^n}f_2(x)dx>0$ can be proved in a similar way.

Let us now consider the case
$$q=q^*=\frac{n}{n-1},\,p>1,\,\,\,\,\text{and}\,\,\,\,\lambda_1\mu_1> 0,\,\mu_1\int_{\mathbb{R}^n}f_1(x)dx>0.$$
As $\mu_1\neq 0,$ from \ref{def2} we have
\begin{equation}\label{3+}\begin{split}
&\int\limits_0^T\int\limits_{\mathbb{R}^n}\left[\frac{\lambda_1}{\mu_1}|u(t, x)|^{p}+|\nabla u(t,x)|^q+\frac{1}{\mu_1}f_1(x)\right]\varphi(t, x)dxdt
\\&\leq \frac{1}{|\mu_1|}\int\limits_0^T\int\limits_{\mathbb{R}^n} |\nabla u(t, x)||\nabla\varphi(t, x)|dxdt+\frac{1}{|\mu_1|}\int\limits_0^T\int\limits_{\mathbb{R}^n}|u(t, x)||\varphi_t(t, x)|dxdt\\&+\frac{1}{|\mu_1|}\int\limits_{\mathbb{R}^{n}}|u_0(x)|\varphi(0, x)dx.\end{split}
\end{equation}
By using the $\varepsilon$-Young inequalities to the right hand side of \eqref{3+} we have
\begin{align*}
\int\limits_0^T\int\limits_{\mathbb{R}^n}|u||\varphi_t|dxdt&\leq \varepsilon_1\int\limits_0^T\int\limits_{\mathbb{R}^n} |u|^p\varphi dxdt +C_1(\varepsilon_1)\int\limits_0^T\int\limits_{\mathbb{R}^n}|\varphi_t|^{\frac{p}{p-1}}{|\varphi|^{-\frac{1}{p-1}}}dxdt,
\end{align*}
and
\begin{align*}
\int\limits_0^T\int\limits_{\mathbb{R}^n}|\nabla u||\nabla\varphi| dxdt&\leq \varepsilon_2\int\limits_0^T\int\limits_{\mathbb{R}^n}|\nabla u|^q\varphi dxdt+C_2(\varepsilon_2)\int\limits_0^T\int\limits_{\mathbb{R}^n}{|\nabla\varphi|^{\frac{q}{q-1}}}{|\varphi|^{-\frac{1}{q-1}}}dxdt,
\end{align*} where $C_1$ and $C_2$ are some positive constants.

Let us choose $\varepsilon_1=\frac{|\lambda_1|}{2},\, \varepsilon_2=\frac{|\mu_1|}{2}$ and inserting the above calculations into \eqref{3+} we have
\begin{equation}\label{3*}\begin{split}
\frac{1}{\mu_1}&\int\limits_0^T\int\limits_{\mathbb{R}^n}f_1(x) \varphi(t, x)dxdt
\leq \frac{1}{|\mu_1|}\int\limits_{\mathbb{R}^{n}}|u_0(x)|\varphi(0, x)dx\\ + C_1&\int\limits_0^T\int\limits_{\mathbb{R}^n}|\varphi_t|^{\frac{p}{p-1}}{|\varphi|^{-\frac{1}{p-1}}}dxdt+ C_2\int\limits_0^T\int\limits_{\mathbb{R}^n}{|\nabla\varphi|^{\frac{q}{q-1}}}{|\varphi|^{-\frac{1}{q-1}}}dxdt.\end{split}
\end{equation}
Let us set the test function as $$\varphi(t,x)=\eta(t)\phi(x),$$
for large enough $R, T$, where $\eta$ and $\phi$ are defined as in \eqref{T2} and \eqref{T1} respectively.

Applying the properties of $\eta$ and $\phi,$ as well as $$|\nabla\phi(x)|\leq C|x|^{-1}\left(\ln R\right)^{-1}$$ and choosing $l\geq\frac{p}{p-1},\,m\geq \frac{q}{q-1}$ we have that
\begin{align*}
&\int\limits_0^T\int\limits_{\mathbb{R}^n}|\varphi_t|^{\frac{p}{p-1}}{|\varphi|^{-\frac{1}{p-1}}}dxdt\leq CT^{1-\frac{p}{p-1}}R^{n},\\& \int\limits_0^T\int\limits_{\mathbb{R}^n}{|\nabla\varphi|^{\frac{q}{q-1}}}{|\varphi|^{-\frac{1}{q-1}}}dxdt\leq CT\left(\ln R\right)^{-\frac{q}{q-1}}.
\end{align*}
Since, $p>1$ and $q=\frac{n}{n-1}$, by inserting the last inequalities into \eqref{3*} and choosing $T=R^j$, we can verify that
\begin{equation}\label{T4*}\begin{split}
&\frac{1}{\mu_1}\int\limits_{\mathbb{R}^n} f_1(x)\phi(x) dx\leq C\left(R^{n-j\frac{p}{p-1}}+(\ln R)^{-n}\right),
\end{split}\end{equation}
Hence choosing $j>\frac{n(p-1)}{p}$ and passing to the limit as $R\to\infty$ in \eqref{T4*} we deduce that \begin{align*}
\frac{1}{\mu_1}\lim\limits_{R\rightarrow +\infty}\int\limits_{\mathbb{R}^n} f_1(x)\phi(x) dx&=\frac{1}{\mu_1}\int\limits_{\mathbb{R}^n} f_1(x) dx\leq 0,\end{align*} which is a contradiction with ${\mu_1}\int\limits_{\mathbb{R}^n} f_1(x) dx>0$.

The case $\mu_2\int_{\mathbb{R}^n}f_2(x)dx>0$ can be proved in a similar way.
The proof is complete.

\section*{Declaration of competing interest}
	The authors declare that there is no conflict of interest.

\section*{Data Availability Statements} The manuscript has no associated data.

\section*{Acknowledgments}
This research has been funded by the Science Committee of the Ministry of Education and Science of the Republic of Kazakhstan (Grant No. AP23483960), by the FWO Odysseus 1 grant G.0H94.18N: Analysis and Partial Differential Equations and by the Methusalem programme of the Ghent University Special Research Fund (BOF) (Grant number 01M01021).


\begin{thebibliography}{AKL92}
\bibitem{Alo} M. Alotaibi, M. Jleli, B. Samet, C. Vetro, First and second critical exponents for an inhomogeneous Schr\"{o}dinger equation with combined nonlinearities. {\it Z. Angew. Math. Phys.}, 73 (2022), 1--17. 

\bibitem{Tao} I. Bejenaru, T. Tao, Sharp well-posedness and ill-posedness results for a quadratic nonlinear Schr\"{o}dinger equation, {\it J. Funct. Anal.}, 233 (2006) 228--259.

\bibitem{Bori} M. Borikhanov, B.T. Torebek, Nonexistence of global solutions for an inhomogeneous pseudo-parabolic equation, {\it Appl. Math. Lett.}, 134 (2022), 108366.

\bibitem{Caze} T. Cazenave, Semilinear Schr\"{o}dinger Equations, Courant Lect. Notes Math., vol. 10, American Mathematical Society, 2003.

\bibitem{Kirane1} A. Z. Fino, I. Dannawi, M. Kirane, Blow-up of solutions for semilinear fractional Schr\"{o}dinger equations. {\it J. Integr. Equ. Appl.}, 30 (2018), 67--80.

\bibitem{Fuj} K. Fujiwara, T. Ozawa, Finite time blowup of solutions to the nonlinear Schr\"{o}dinger equation without gauge invariance, {\it J. Math. Phys.}, 57 (2016) 082103.

\bibitem{Fuj1} K. Fujiwara, T. Ozawa, Lifespan of strong solutions to the periodic nonlinear Schr\"{o}dinger equation without gauge invariance. {\it J. Evol. Equ.}, 17 (2017), 1023--1030.

\bibitem{Ikeda1} M. Ikeda, T. Inui, Small data blow-up of $L^2$ or $H^1$-solution for the semilinear Schr\"{o}dinger equation without gauge invariance. {\it J. Evol. Equ.}, 15 (2015), 571--581.

\bibitem{Ikeda2} M. Ikeda, T. Inui, Some non-existence results for the semilinear Schr\"{o}dinger equation without gauge invariance. {\it J. Math. Anal. Appl.}, 425 (2015), 758--773.

\bibitem{Ikeda} M. Ikeda, Y. Wakasugi, Small-data blow-up of $L^2$-solution for the nonlinear Schr\"{o}dinger equation without gauge invariance. {\it Differential Integral Equations}, 26 (2013), 1275--1285.

\bibitem{Jleli} M. Jleli, B. Samet, On the critical exponent for nonlinear Schr\"{o}dinger equations without gauge invariance in exterior domains. {\it J. Math. Anal. Appl.}, 469 (2019), 188--201.

\bibitem{Souplet} M. Jleli, B. Samet, P. Souplet, Discontinuous critical Fujita exponents for the heat equation with combined nonlinearities, {\it Proc. Am. Math. Soc.}, 148 (2020), 2579--2593.

\bibitem{Samet} M. Jleli, B. Samet, C. Vetro, A blow-up result for a nonlinear wave equation on manifolds: the critical case, {\it Appl. Anal.}, 102 (2023) 1463--1472.

\bibitem{Kirane} M. Kirane, A. Nabti, Life span of solutions to a nonlocal in time nonlinear fractional Schr\"{o}dinger equation. {\it Z. Angew. Math. Phys.}, 66 (2015), 1473--1482.

\bibitem{Pohoz} E. Mitidieri, S.I. Pohozaev, A priori estimates and the absence of solutions of nonlinear partial differential equations and inequalities, {\it Tr. Mat. Inst. Steklova}, 234 (2001) 1--384 (in Russian).

\bibitem{Oh} T. Oh, A blowup result for the periodic NLS without gauge invariance, {\it C. R. Acad. Sci. Paris}, 350 (2012), 389--392.

\bibitem{Torebek} B. T. Torebek, Critical exponents for the $p$-Laplace heat equations with combined nonlinearities, {\it J. Evol. Equ.}, 23 (2023), 1--15.

\bibitem{Zhang} Q. S. Zhang, A new critical phenomenon for semilinear parabolic problems, {\it J. Math. Anal. Appl.}, 219  (1998), 125--139.
\bibitem{Zhang1} Q. S. Zhang, Blow-up results for nonlinear parabolic equations on manifolds, {\it Duke Math. J.}, 97 (1999), 515--539.
\end{thebibliography}
\end{document}